\documentclass[11pt]{article}
\usepackage{amsfonts}
\usepackage{amscd}
\usepackage{amsmath}
\usepackage{epsfig}

\newtheorem{theorem}{Theorem}
\newtheorem{proposition}{Proposition}

\newtheorem{remark}{Remark}
\newtheorem{lemma}{Lemma}

\def\cy{s}
\def\ko{\zeta}

\def\be{\begin{equation}}
\def\ee{\end{equation}}
\def\ben{\begin{displaymath}}
\def\een{\end{displaymath}}
\def\baa{\begin{eqnarray}}
\def\eaa{\end{eqnarray}}

\def\ba{\begin{array}}
\def\ea{\end{array}}
\def\CP1{{\bf CP}^1}

\def\a{\alpha}

\def\b{\beta}

\def\l{\lambda}
\def\r{\rho}

\def\P{\wp}
\def\Fcal{{\cal F}}

\def\f{\frac}
\def\la{\label}

\def\L{{\cal X}}

\def\p{\partial}
\def\B{{\bf B}}

\def\Be{{\bf w}}

\def\c{{\cal C}}

\def\P{P}

\def\cdiff{{\cal C}}

\def\Abel{{\cal A}}
\def\hd{v}
\def\cdiff{{\cal C}}
\def\Wcal{{\cal W}}

\def\g{{\bf m}}

\begin{document}

\title{Compact polyhedral surfaces of an arbitrary genus and determinants of Laplacians}
\author{Alexey Kokotov
\footnote{Department of Mathematics and Statistics, Concordia University,
1455 de Maisonneuve Blvd. West, Montreal, Quebec, H3G 1M8 Canada,
{\bf E-mail: alexey@mathstat.concordia.ca}}} \maketitle

 \vskip0.5cm {\bf Abstract.} Compact polyhedral surfaces (or, equivalently,
 compact Riemann surfaces with conformal flat conical metrics) of an arbitrary
 genus are considered. After giving a short self-contained survey of their
 basic spectral properties, we study the zeta-regularized determinant of the
 Laplacian as a functional on the moduli space of these surfaces. An explicit
 formula for this determinant is obtained. \vskip0.5cm

\section{Introduction}

There are several well-known ways to introduce a compact Riemann surface, e.g., via algebraic equations
or by means of some uniformization theorem, where the surface is introduced
as the quotient of the upper half-plane over the action of a Fuchsian group.
In this paper we consider a  less popular approach which is at the same
time, perhaps, the most elementary:  one can simply consider the boundary of
a connected (but, generally, not simply connected) polyhedron in three
dimensional Euclidean space. This is a polyhedral surface which carries the
structure of a complex manifold (the corresponding system of holomorphic
local parameters is obvious for all points except the vertices; near a vertex
one should introduce the local parameter $\zeta=z^{2\pi/\alpha}$, where
$\alpha$ is the sum of the angles adjacent to the vertex). In this way the
Riemann surface arises together with a conformal metric; this metric is flat
and has conical singularities at the vertices. Instead of a polyhedron one
can also start from some abstract simplicial complex, thinking of a
polyhedral surface as glued from plane triangles.

The present paper is devoted to the spectral theory of the Laplacian on such
surfaces. The main goal is to study the determinant of the Laplacian (acting
in the trivial line bundle over the surface) as a functional on the space of
Riemann surfaces with conformal flat conical metrics (polyhedral surfaces).
The similar question for {\it smooth} conformal metrics and arbitrary
holomorphic bundles was very popular in the eighties and early nineties being
motivated by string theory. The determinants of Laplacians in flat
singular metrics  are much less studied: among the very few appropriate
references we mention \cite{DP}, where the determinant of the Laplacian in
a conical metric was defined via some special regularization of the diverging
Liouville integral and the question about the relation of such a definition
with the spectrum of the Laplacian remained open, and two papers \cite{King},
\cite{AuSal} dealing with flat conical metrics on the Riemann sphere.

In  \cite{Leipzig} (see also \cite{Leipzig1}) the determinant of the Laplacian was studied as a
functional
$${\cal H}_g(k_1, \dots, k_M)\ni (\L, \omega)\mapsto {\rm det}\,\Delta^{|\omega|^2}$$
on the space ${\cal H}_g(k_1, \dots, k_M)$ of equivalence classes of pairs
$(\L, \omega)$, where $\L$ is a compact Riemann surface of genus $g$ and
$\omega$ is a holomorphic one-form (an Abelian differential) with $M$ zeros
of multiplicities $k_1, \dots, k_M$. Here ${\rm det}\,\Delta^{|\omega|^2}$
stands for the determinant of the Laplacian in the flat metric $|\omega|^2$
having conical singularities at the zeros of $\omega$. The flat conical
metric $|\omega|^2$ considered in \cite{Leipzig} is very special: the divisor
of the conical points of this metric is not arbitrary (it should be the
canonical one, i. e. coincide with the divisor of a holomorphic one-form) and
the conical angles at the conical points are integer multiples of $2\pi$.
Later in \cite{Klochko} this restrictive condition has been eliminated in the
case of polyhedral surfaces of genus one.

In the present paper we generalize the results of \cite{Leipzig} and
\cite{Klochko} to the case of polyhedral surfaces of an arbitrary genus.
Moreover, we give a short and self-contained survey of some basic facts from
the spectral theory of the Laplacian on flat surfaces with conical points. In
particular, we discuss the theory of self-adjoint extensions of this
Laplacian and study the asymptotics of the corresponding heat kernel.

\section{Flat conical metrics on surfaces}
Following \cite{Tro} and \cite{Klochko}, we discuss here flat conical metrics
on compact Riemann surfaces of an arbitrary genus.
\subsection{Troyanov's theorem}\index{Troyanov theorem on flat metrics with conical singularities}
Let $\sum_{k=1}^Nb_kP_k$ be a (generalized, i.e., the coefficients $b_k$ are
not necessary integers) divisor on a compact Riemann surface $\L$ of genus
$g$. Let also $\sum_{k=1}^Nb_k=2g-2$. Then, according to Troyanov's theorem
(see \cite{Tro}),  there exists a (unique up to a rescaling) conformal (i. e.
giving rise to a complex structure which coincides with that of $\L$) flat
metric $\g$ on $\L$ which is smooth in $\L\setminus\{P_1, \dots, P_N\}$ and
has simple singularities of order $b_k$ at $P_k$. The latter means that in a
vicinity of $P_k$ the metric $\g$ can be represented in the form
\begin{equation}\label{met}\g=e^{u(z, \bar z)}|z|^{2b_k}|dz|^2,\end{equation}
where $z$ is a conformal coordinate and $u$ is a smooth real-valued function.
In particular, if $\beta_k>-1$ the point $P_k$ is conical with conical angle
$\beta_k=2\pi (b_k+1)$. Here we construct the metric $\g$ explicitly, giving
an effective proof of Troyanov's theorem (cf. \cite{Klochko}).

Fix a canonical basis of cycles on $\L$ (we assume that $g\geq 1$, the case
$g=0$ is trivial) and let $E(P, Q)$ be the prime-form (see \cite{Fay73}).
Then for any divisor ${\cal D}=r_1Q_1+\dots r_mQ_M-s_1R_1-\dots -s_NR_N$ of
degree zero on $\L$ (here the coefficients $r_k, s_k$ are positive integers)
the meromorphic differential
$$\omega_{{\cal D}}=d_z\log\frac{\prod_{k=1}^ME^{r_k}(z,Q_k)}{\prod_{k=1}^NE^{s_k}(z, R_k)}$$
is holomorphic outside ${\cal D}$ and has first order poles at the points
of ${\cal D}$ with residues $r_k$ at $Q_k$ and $-s_k$ at $R_k$. Since the
prime-form is single-valued along the ${\bf a}$-cycles, all ${\bf
a}$-periods of the differential $\omega_{\cal D}$ vanish.

Let $\{v_\alpha\}_{\alpha=1}^g$ be the basis of holomorphic normalized
differentials and ${\mathbb B}$ the corresponding matrix of ${\bf
b}$-periods. Then all ${\bf a}$- and ${\bf b}$-periods of the meromorphic
differential
$$\Omega_{\cal D}=\omega_{\cal D}-2\pi i\sum_{\a, \b=1}^g ((\Im{\mathbb B})^{-1})_{\a\b}\Im\left( \int_{s_1R_1+\dots s_NR_N}^{r_1Q_1+\dots r_M Q_M}v_\beta\right)v_\alpha$$
are purely imaginary (see \cite{Fay73}, p. 4).

Obviously, the differentials $\omega_{\cal D}$ and $\Omega_{\cal D}$ have the
same structure of poles: their difference is a holomorphic $1$-form.

Choose a base-point $P_0$ on $\L$ and introduce the following quantity
$${\cal F}_{\cal D}(P)=\exp\int_{P_0}^P\Omega_{\cal D}.$$
Clearly, ${\cal F}_{\cal D}$ is a meromorphic section of some {\it unitary}
flat line bundle over $\L$, the divisor of this section coincides with ${\cal
D}$.

Now we are ready to construct the metric $\g$. Choose any holomorphic
differential $w$ on $\L$ with, say, only simple zeros $S_1, \dots, S_{2g-2}$.
Then one can set $\g=|u|^2$, where
\begin{equation}\label{reshen}u(P)=w(P)\Fcal_{(2g-2)S_0-S_1-\dots S_{2g-2}}(P)
\prod_{k=1}^N\left[\Fcal_{P_k-S_0}(P)\right]^{b_k}\end{equation} and $S_0$ is
an arbitrary point.

Notice that in the case $g=1$ the second factor in (\ref{reshen}) is absent and
the remaining part is nonsingular at the point $S_0$.

\subsection{Distinguished local parameter} \index{distinguished
local parameter} In a vicinity of a conical point the flat metric (\ref{met})
takes the form
$$\g=|g(z)|^2|z|^{2b}|dz|^2$$
with some holomorphic function $g$ such that $g(0)\neq0$. It is easy to show
(see, e. g., \cite{Tro},  Proposition 2) that there exists a holomorphic
change of variable $z=z(x)$ such that in the local parameter $x$
$$\g=|x|^{2b}|dx|^2\,.$$
We shall call the parameter $x$ (unique up to a constant factor $c$, $|c|=1$)
{\it distinguished}. In case $b>-1$ the existence of the  distinguished
parameter means that in a vicinity of a conical point the surface $\L$ is
isometric to the standard cone with conical angle $\beta=2\pi(b+1)$.

\subsection {Euclidean polyhedral surfaces.} In \cite{Tro} it is proved that
 any compact Riemann surface with flat conformal conical metric admits
 a proper triangulation (i. e. each conical point is a vertex of some triangle of the triangulation).
This means that any compact Riemann surface with a flat conical metric is a
{\it Euclidean polyhedral surface} (see \cite{Bobenko}) i. e. can be glued
from Euclidean triangles. On the other hand as it is explained in
\cite{Bobenko} any compact Euclidean oriented polyhedral surface gives rise
to a Riemann surface with a flat conical metric. Therefore, from now on we do
not discern compact Euclidean polyhedral surfaces and Riemann surfaces with
flat conical metrics.

\section{Laplacians on polyhedral surfaces. Basic facts}
Without claiming originality we give here a short  self-contained survey of some
basic facts from the spectral theory of Laplacian on compact polyhedral
surfaces. We start with recalling the (slightly modified) Carslaw
construction (1909) of the heat kernel on a cone, then we describe the set of
self-adjoint extensions of a conical Laplacian (these results are complementary
to Kondratjev's study (\cite{Kondr}) of elliptic equations on conical
manifolds and are well-known, being in the folklore since the sixties
of the last century; their
generalization to the case of Laplacians acting on $p$-forms can be found in
\cite{Mooers}). Finally, we establish the precise heat asymptotics for the
Friedrichs extension of the Laplacian on a compact polyhedral surface. It
should be noted that more general results on the heat asymptotics for
Laplacians acting on $p$-forms on piecewise flat pseudomanifolds can be found
in \cite{Cheeger}.

\subsection{The heat kernel on the infinite cone}

 We start from the standard heat kernel
  \begin{equation}\label{h1}H_{2\pi}( x,  y; t)=\frac{1}{4\pi
 t}\exp\{-(x-y)\cdot(x-y)/4t\}\end{equation}
 in the space ${\mathbb R}^2$ which we consider as the cone with conical angle $2\pi$.
  Introducing the polar coordinates $(r, \theta)$ and $(\rho,
 \psi)$ in the $x$ and $y$-planes, one can rewrite (\ref{h1}) as the contour integral
$$H_{2\pi}( x,  y;
t)=$$
\begin{equation}\label{h2} \frac{1}{16\pi^2it}\exp\{-(r^2+\rho^2)/4t\}\int_{C_{\theta,
\psi}}\exp\{r\rho\cos(\alpha-\theta)/2t\}\cot
\frac{\alpha-\psi}{2}\,d\alpha,\end{equation} where $C_{\theta, \psi}$
denotes the union of a small positively oriented circle centered at
$\alpha=\psi$ and the two vertical lines, $l_1=(\theta-\pi-i\infty,
\theta-\pi+i\infty)$ and $l_2=(\theta+\pi+i\infty, \theta+\pi-i\infty)$,
having mutually opposite orientations.

To prove (\ref{h2}) one has to notice that

1) $\Re \cos(\alpha-\theta)<0$ in vicinities of the lines $l_1$ and $l_2$
and, therefore, the integrals over these lines converge.

2) The integrals over the lines cancel due to the $2\pi$-periodicity of the
integrand and the remaining integral over the circle coincides with
(\ref{h1}) due to the Cauchy Theorem.

Observe that one can deform the contour $C_{ \theta, \psi}$ into the union,
$A_{\theta}$, of two contours lying in the open domains
$\{\theta-\pi<\Re\alpha<\theta+\pi\,,\, \Im \alpha>0\}$ and
$\{\theta-\pi<\Re\alpha<\theta+\pi\,,\,  \Im \alpha<0\}$ respectively, the
first contour goes from $\theta+\pi+i\infty$ to $\theta-\pi+i\infty$, the
second one goes from $\theta-\pi-i\infty$ to $\theta+\pi-i\infty$. This leads
to the following representation for the heat kernel $H_{2\pi}$:
 $$H_{2\pi}( x,  y;
t)=$$
\begin{equation}\label{h3}\frac{1}{16\pi^2it}\exp\{-(r^2+\rho^2)/4t\}\int_{A_
\theta}\exp\{r\rho\cos(\alpha-\theta)/2t\}\cot
\frac{\alpha-\psi}{2}\,d\alpha.\end{equation}

The latter representation admits a natural generalization to the case of the
cone $C_\beta$ with conical angle $\beta$, $0<\beta<+\infty$. Notice here
that in case $0<\beta\leq 2\pi$ the cone $C_\beta$ is isometric to the
surface $ z_3=\sqrt{(\frac{4\pi^2}{\beta^2}-1)(z_1^2+z_2^2)}$.

Namely, introducing the polar coordinates on $C_\beta$, we see that the
following expression represents the heat kernel on $C_\beta$:
$$H_{\beta}(r,
\theta, \rho, \psi; t)=$$
\begin{equation}\label{h4}
\frac{1}{8\pi\beta it}\exp\{-(r^2+\rho^2)/4t\}\int_{A_
\theta}\exp\{r\rho\cos(\alpha-\theta)/2t\}\cot
\frac{\pi(\alpha-\psi)}{\beta}\,d\alpha\,.\end{equation}

Clearly, expression (\ref{h4}) is symmetric with respect to $(r, \theta)$ and
$(\rho, \psi)$ and is $\beta$-periodic with respect to the angle variables
$\theta, \psi$. Moreover, it satisfies the heat equation on $C_\beta$.
Therefore,  to verify that $H_{\beta}$ is in fact the heat kernel on
$C_\beta$ it remains to show that $H_\beta(\cdot, y, t)\longrightarrow
\delta(\cdot-y)$ as $t\to 0+$. To this end deform the contour $A_\psi$ into
the union of the lines $l_1$ and $l_2$ and (possibly many) small circles
centered at the poles of $\cot \frac{\pi(\cdot-\psi)}{\beta}$ in the strip
$\theta-\pi<\Re\alpha<\theta+\pi$. The integrals over all the components of
this union except the circle centered at $\alpha=\psi$ vanish in the limit as
$t\to 0+$, whereas the integral over the latter circle coincides with
$H_{2\pi}$. \subsubsection{The heat asymptotics near the vertex}
\begin{proposition}\label{ugol}
Let $R>0$ and $C_{\beta}(R)=\{x\in C_{\beta}: {\rm dist}(x, {\cal O})< R\}$.
Let also $dx$ denote the area element on $C_\beta$. Then for some
$\epsilon>0$
\begin{equation}\label{as1}\int_{C_\beta(R)}H_{\beta}(x, x; t)\,dx=\frac{1}{4\pi
t}{\rm Area}(C_{\beta}(R))+
\frac{1}{12}\left(\frac{2\pi}{\beta}-\frac{\beta}{2\pi}
\right)+O(e^{-\epsilon/t})\end{equation} as $t\to 0+$.
\end{proposition}

{\bf Proof}\,(cf. \cite{Fursaev}, p. 1433). Make in (\ref{h4}) the change of
variable $\gamma=\alpha-\psi$ and
 deform the contour
$A_{\theta-\psi}$ into the contour $\Gamma^-_{\theta-\psi}\cup
\Gamma^{+}_{\theta-\psi}\cup \{ |\gamma|=\delta\}$, where the oriented curve
$\Gamma^{-}_{\theta-\psi}$ goes from $\theta-\psi-\pi-i\infty$ to
$\theta-\psi-\pi+i\infty$ and intersects the real axis at $\gamma=-\delta$,
the oriented curve $\Gamma^+_{\theta-\psi}$ goes from
$\theta-\psi+\pi+i\infty$ to $\theta-\psi+\pi-i\infty$ and intersects the
real axis at $\gamma=\delta$, the circle $\{|\gamma|=\delta\}$ is positively
oriented and $\delta$ is a small positive number. Calculating the integral
over the circle $\{|\gamma|=\delta\}$  via the Cauchy Theorem, we get
$$H_\beta(x, y; t)-H_{2\pi}(x, y; t)=$$
\begin{equation}\label{h6}
\frac{1}{8\pi\beta
it}\exp\{-(r^2+\rho^2)/4t\}\int_{\Gamma^-_{\theta-\psi}\cup\Gamma^+_{\theta-\psi}}\exp\{r\rho\cos(\gamma+\psi-\theta)/2t\}\cot
\frac{\pi\gamma}{\beta}\,d\gamma
\end{equation}
and
$$\int_{C_\beta(R)}\left(H_\beta(x, x; t)-\frac{1}{4\pi t}\right)dx=$$
\begin{equation}\label{h7}
\frac{1}{8\pi i t}\int_0^R\,dr\,
r\int_{\Gamma_0^-\cup\Gamma_0^+}\exp\{-\frac{r^2\sin^2(\gamma/2)}{t}\}\cot\frac{\pi\gamma}{\beta}\,d\gamma\,.
\end{equation}
 The  integration over $r$ can be done
explicitly  and the right hand side of (\ref{h7}) reduces to
\begin{equation}\label{h8}\frac{1}{16\pi i}
\int_{\Gamma_0^-\cup\Gamma_0^+}\frac{\cot(\frac{\pi\gamma}{\beta})}{\sin^2(\gamma/2)}\,d\gamma+O(e^{-\epsilon/t}).
\end{equation}
(One can assume that $\Re\sin^2(\gamma/2)$ is positive and separated from
zero when $\gamma\in\Gamma_0^-\cup\Gamma_0^+$.)
 The contour of integration in (\ref{h8}) can be changed for a
negatively oriented circle centered at $\gamma=0$. Since ${\rm Res}
(\frac{\cot(\frac{\pi\gamma}{\beta})}{\sin^2(\gamma/2)}\,,\,\gamma=0)=\frac{2}{3}(\frac{\beta}{2\pi}-\frac{2\pi}{\beta})$,
we arrive at (\ref{as1}).

\begin{remark} {\rm The Laplacian $\Delta$ corresponding to the flat conical metric $(d\r)^2+r^2(d\theta)^2,  0\leq\theta\leq
\beta$ on $C_\beta$ with domain $C^{\infty}_{0}(C_\beta\setminus {\cal O})$
has infinitely many self-adjoint extensions. Analyzing the asymptotics of
(\ref{h4}) near the vertex ${\cal O}$, one can show that for any $y\in
C_\beta, t>0$ the function $H_\beta(\cdot, y; t)$ belongs to the domain of
the Friedrichs extension $\Delta_F$ of $\Delta$ and does not belong to the
domain of any other extension. Moreover, using a Hankel transform, it is
possible to get an explicit spectral representation of $\Delta_F$ (this
operator has an absolutely continuous spectrum of infinite multiplicity) and to
show that the Schwartz kernel of the operator $e^{t\Delta_F}$ coincides with
$H_{\beta}(\cdot, \cdot; t)$ (see, e. g., \cite{Taylor} formula (8.8.30)
together with \cite{Carslaw}, p. 370.)
 }\end{remark}

\subsection{Heat asymptotics for compact polyhedral surfaces}
\subsubsection{Self-adjoint extensions of a conical Laplacian}
Let $\L$ be a compact polyhedral surface with vertices (conical points) $P_1,
\dots, P_N$. The Laplacian $\Delta$ corresponding to the natural flat conical
metric on $\L$ with domain $C^{\infty}_0(\L\setminus\{P_1, \dots, P_N\})$ (we
remind the reader that the Riemannian manifold $\L$ is smooth everywhere
except the vertices) is not essentially self-adjoint and  one has to fix one
of its self-adjoint extensions. We are to discuss now the choice of a
self-adjoint extension.

This choice is defined by the prescription of some particular asymptotical
behavior  near the conical points to functions from the domain of the
Laplacian; it is sufficient to consider a surface with  only one conical
point $P$ of the conical angle $\beta$. More precisely, assume that $\L$ is
smooth everywhere except the point $P$ and that some vicinity of $P$ is
isometric to a vicinity of the vertex ${\cal O}$ of the standard cone
$C_{\beta}$ (of course, now the metric on $\L$ no more can be flat everywhere
in $\L\setminus P$ unless the genus $g$ of $\L$ is greater than one and
$\beta=2\pi(2g-1)$).

For $k\in {\mathbb N}_0$ introduce the functions  $V_\pm^k$ on $C_\beta$ by
$$V^k_\pm(r, \theta)=r^{{\pm}\frac{2\pi k}{\beta}}\exp\{i\frac{2\pi
k\theta}{\beta} \};\ \ k>0\,,$$
$$V_+^0=1,\ \  V_-^0=\log r\,.$$
Clearly, these functions are formal solutions to the homogeneous problem
$\Delta u=0$ on $C_\beta$. Notice that the functions $V_-^k$ grow near the
vertex but are still square integrable in its vicinity if
$k<\frac{\beta}{2\pi}$.

Let ${\cal D}_{\rm min}$ denote the graph closure of $C^\infty_0(\L\setminus
P)$, i.e.,
$$U\in {\cal D}_{{\rm min}} \Leftrightarrow \exists u_m\in C^\infty_0(\L\setminus P), W\in L_2(\L):
u_m\rightarrow U \ {\rm and}\ \Delta u_m \rightarrow W\ \ {\rm in}\ \
L_2(\L).$$

Define the space $H_\delta^2(C_\beta)$ as the closure of
$C^\infty_0(C_\beta\setminus{\cal O})$ with respect to the norm
$$||u; H_\delta^2(C_\beta)||^2=\sum_{|{\bf \alpha}|\leq
2}\int_{C_\beta}r^{2(\delta-2+|{\bf \alpha}|)}|D^{\bf \alpha}_xu(x)|^2dx.$$
Then for any $\delta\in {\mathbb R}$ such that $\delta-1\neq \frac{2\pi
k}{\beta}, k\in {\mathbb Z}$ one has the a priori estimate
\begin{equation}\label{S1}||u; H^2_\delta(C_\beta)||\leq c||\Delta u;
H^0_\delta(C_\beta)||\end{equation} for any $u\in
C^\infty_0(C_\beta\setminus{\cal O})$ and some constant $c$ being independent
of $u$ (see, e.g., \cite{NazPla}, Chapter 2).

It follows from Sobolev's imbedding theorem that for functions from $u\in
H_\delta^2(C_\beta)$ one has the point-wise estimate
\begin{equation}\label{S2} r^{\delta-1}|u(r, \theta)|\leq c||u;
H_\delta^2(C_\beta)||.\end{equation}

Applying estimates (\ref{S1}) and (\ref{S2}) with $\delta=0$, we see that
functions $u$ from ${\cal D}_{\rm min}$ must obey the asymptotics $u(r,
\theta) =O(r)$ as $r\to 0$.

Now the description of the set of all self-adjoint extensions of $\Delta$
looks as follows. Let $\chi$ be a smooth function on $\L$ which is equal to
$1$ near the vertex $P$ and such that in a vicinity of the support of $\chi$
$\L$ is isometric to $C_\beta$. Denote by  ${\frak M}$ the linear subspace of
$L_2(\L)$ spanned by the functions $\chi V_\pm^k$ with $0\leq k<
\frac{\beta}{2\pi}$. The dimension, $2d$, of ${\frak M}$ is even. To get a
self-adjoint extension of $\Delta$ one chooses a subspace ${\frak N}$ of
${\frak M}$ of dimension $d$ such that
 $$(\Delta u, v)_{L_2(\L)}-(u, \Delta
v)_{L_2(\L)}=\lim_{\epsilon\to 0+}\oint_{r=\epsilon}\left(u\frac{\partial
v}{\partial r}-v\frac{\partial u}{\partial r}\right)=0$$
 for any $u, v\in {\frak N}$. To any such subspace ${\frak N}$
there corresponds a self-adjoint extension $\Delta_{\frak N}$ of $\Delta$
with domain ${\frak N}+{\cal D}_{{\rm min}}$.

The extension corresponding to the subspace ${\frak N}$ spanned by the
functions $\chi V_+^k$, $0\leq k<\frac{\beta}{2\pi}$ coincides with the
Friedrichs extension of $\Delta$. The functions from the domain of the
Friedrichs extension are bounded near the vertex.

From now on we denote by $\Delta$ the Friedrichs extension of the Laplacian
on the polyhedral surface $\L$; other extensions will not be considered here.

\subsubsection{Heat asymptotics}\index{heat asymptotics for Euclidean polyhedral surfaces}
\begin{theorem}\label{Mainsp}
Let $\L$ be a compact polyhedral surface with vertices $P_1, \dots, \P_N$ of
conical angles $\beta_1, \dots, \beta_N$. Let $\Delta$ be the Friedrichs
extension of the Laplacian defined on functions from $C^\infty_0(\L\setminus
\{P_1, \dots, P_N\})$. Then

\begin{enumerate}\item The spectrum of the operator $\Delta$ is discrete,
all the eigenvalues of $\Delta$ have finite multiplicity.
\item Let ${\cal H}(x, y; t)$ be the heat kernel for $\Delta$. Then for some $\epsilon>0$
\begin{equation}\label{Mainasym}
{\rm Tr}\,e^{t\Delta}=\int_\L {\cal H}(x, x; t)\,dx=\frac{{\rm
Area}(\L)}{4\pi
t}+\frac{1}{12}\sum_{k=1}^N\left\{\frac{2\pi}{\beta_k}-\frac{\beta_k}{2\pi}
\right\}+O(e^{-\epsilon/t}),
\end{equation}
as $t\to 0+$.
\item The counting function, $N(\lambda)$, of the spectrum of $\Delta$ obeys
the asymptotics $N(\lambda)=O(\lambda)$ as $\lambda\to +\infty$.
\end{enumerate}
\end{theorem}
{\bf Proof.} 1) The proof of the first statement is a standard exercise (cf.
\cite{King}). We indicate only the main idea leaving the details to the
reader.
 Introduce the closure, $H^1(\L)$, of $C^\infty_0(\L\setminus\{P_1, \dots, P_N\}$ with respect to the norm $|||u|||=||u;
L_2||+||\nabla u; L_2||$. It is sufficient to prove that any  bounded set $S$
in $H^1(\L)$ is precompact in the $L_2$-topology (this will imply the compactness
of the self-adjoint operator $(I-\Delta)^{-1}$). Moreover, one can assume
that the supports of functions from $S$ belong to a small ball $B$ centered
at a conical point $P$. Now to prove the precompactness of $S$ it is
sufficient to make use of the expansion with respect to eigenfunctions of the
Dirichlet problem in $B$ and the diagonal process.

2)Let $\L=\cup_{j=0}^NK_j$, where $K_j$, $j=1, \dots, N$ is a neighborhood of
the conical point $P_j$ which is isometric to $C_{\beta_j}(R)$ with some
$R>0$, and $K_0=\L\setminus\cup_{j=1}^NK_j$.

Let also $K^{\epsilon_1}_j\supset K_j$ and $K^{\epsilon_1}_j$ be isometric to
$C_{\beta_j}(R+\epsilon_1)$ with some $\epsilon_1>0$ and $j=1, \dots, N$.

Fixing $t>0$ and $x, y\in K_j$ with $j>0$, one has
\begin{equation}\label{chasti}\int_0^t\,ds\int_{K_j^{\epsilon_1}}\left(
\psi\{\Delta_z-\partial_s\}\phi-\phi\{\Delta_z+\partial_s\}\psi\right)\,dz=\end{equation}
\begin{equation}
\int_0^tds\int_{\partial K_j^{\epsilon_1}}\left(\phi\frac{\partial
\psi}{\partial n}-\psi\frac{\partial \phi}{\partial
n}\right)dl(z)-\int_{K_j^{\epsilon_1}}\left(\phi(z, t)\psi(z,t)-\phi(z,
0)\psi(z, 0)\right)\,dz\end{equation}
 with $\phi(z,
t)={\cal H}(z, y; t)-H_{\beta_j}(z, y; t)$ and $\psi(z, t)=H_{\beta_j}(z, x;
t-s)$. (Here it is important that we are working with the heat kernel of the
Friedrichs extension of the Laplacian, for other extensions the heat kernel
has growing terms in the asymptotics near the vertex and the right hand side
of (\ref{chasti}) gets extra terms.) Therefore,
$$H_{\beta_j}(x, y; t)-{\cal H}(x, y; t)=$$
$$\int_0^tds\int_{\partial K_j^{\epsilon_1}}\left({\cal H}(y, z; s)\frac{\partial H_{\beta_j}(x, z; t-s)}{\partial n(z)}
-H_{\beta_j}(z, x; t-s)\frac{{\partial \cal H}(z, y; s)}{\partial
n(z)}\right)\,dl(z)$$$$=O(e^{-\epsilon_2/t})$$ with some $\epsilon_2>0$ as
$t\to 0+$ uniformly with respect to $x, y\in K_j$. This implies that
\begin{equation}\label{assy1}\int_{K_j}{\cal H}(x, x; t)dx=\int_{K_j} H_{\beta_j}(x, x;
t)dx+O(e^{-\epsilon_2/t}).\end{equation} Since the metric on $\L$ is flat in
a vicinity of  $K_0$, one has the asymptotics
$$\int_{K_0}{\cal H}(x, x; t)dx=\frac{{\rm Area}(K_0)}{4\pi t}+O(e^{-\epsilon_3/t})$$
with some $\epsilon_3>0$ (cf. \cite{McKeanSinger}). Now (\ref{Mainasym})
follows from (\ref{as1}).

3) The third statement of the theorem follows from the second one due to the
standard Tauberian arguments.

\section{Determinant of the Laplacian: Analytic surgery and Polyakov-type
formulas}\index{zeta-regularized determinant of Laplacian}

Theorem \ref{Mainsp} opens a way to define the determinant, ${\rm
det}\,\Delta$, of the Laplacian on a compact polyhedral surface via the
standard Ray-Singer regularization. Namely introduce the operator
$\zeta$-function
\begin{equation}\label{oper}\zeta_{\Delta}(s)=\sum_{\l_k>0}\frac{1}{\lambda_k^s},\end{equation} where the
summation goes over all strictly positive eigenvalues $\l_k$ of the operator
$-\Delta$ (counting  multiplicities). Due to the third statement of Theorem
\ref{Mainsp}, the function $\zeta_\Delta$ is holomorphic in the half-plane
$\{\Re s>1\}$. Moreover, due to the equality
\begin{equation}\label{zeta}\zeta_\Delta(s)=\frac{1}{\Gamma(s)}\int_0^\infty\left\{{\rm
Tr}\,e^{t\Delta }-1\right\}t^{s-1}\,dt\end{equation} and the asymptotics
(\ref{Mainasym}), one has the equality
\begin{equation}\label{predst}
\zeta_\Delta(s)=\frac{1}{\Gamma(s)}\left\{ \frac{{\rm
Area}\,(\L)}{4\pi(s-1)}+\left[\frac{1}{12}\sum_{k=1}^N\left\{\frac{2\pi}{\beta_k}-\frac{\beta_k}{2\pi}
\right\}-1\right]\frac{1}{s}+e(s)
 \right\},
\end{equation}
where $e(s)$ is an entire function. Thus, $\zeta_\Delta$ is regular at $s=0$
and one can define the $\zeta$-regularized determinant of the Laplacian via
usual $\zeta$-regularization (cf. \cite{Ray}):
 \begin{equation}\label{defdet}{\rm
det}\Delta:=\exp\{-\zeta'_\Delta(0)\}\,.\end{equation}
 Moreover,
(\ref{predst}) and the relation $\sum_{k=1}^Nb_k=2g-2$;
$b_k=\frac{\beta_k}{2\pi}-1$ yield
\begin{equation}\label{euler}
\zeta_\Delta(0)=\frac{1}{12}\sum_{k=1}^N\left\{\frac{2\pi}{\beta_k}-\frac{\beta_k}{2\pi}
\right\}-1=\left(\frac{\chi(\L)}{6}-1\right)+\frac{1}{12}\sum_{k=1}^N\left\{\frac{2\pi}{\beta_k}+\frac{\beta_k}{2\pi}
-2\right\},
\end{equation}
where $\chi(\L)=2-2g$ is the Euler characteristics of $\L$.

It should be noted that the term $\frac{\chi(\L)}{6}-1$ at the right hand
side of (\ref{euler}) coincides with the value at zero of the operator
$\zeta$-function of the Laplacian corresponding to an arbitrary {\it smooth}
metric on $\L$ (see, e. g., \cite{Sarnak}, p. 155).

Let $\g$ and $\tilde{\g}=\kappa\g$, $\kappa>0$ be two homothetic flat metrics
with the same conical points with conical angles $\beta_1, \dots, \beta_N$.
Then  (\ref{oper}), (\ref{defdet}) and (\ref{euler}) imply the following {\it
rescaling property} of the conical Laplacian:
\begin{equation}\label{rescaling}
{\rm
det}\Delta^{\tilde{\g}}=\kappa^{-\left(\frac{\chi(\L)}{6}-1\right)-\frac{1}{12}\sum_{k=1}^N\left\{\frac{2\pi}{\beta_k}+\frac{\beta_k}{2\pi}
-2\right\}}{\rm det}\,\Delta^{\g}
\end{equation}

\subsection{Analytic surgery}\label{AASS}\index{Analytic surgery}

Let $\g$ be an arbitrary smooth metric on $\L$ and denote by $\Delta^\g$ the
corresponding Laplacian. Consider $N$ nonoverlapping connected and simply
connected domains $D_1, \dots, D_N\subset \L$ bounded by closed curves
$\gamma_1, \dots, \gamma_N$ and introduce also the domain
$\Sigma=\L\setminus\cup_{k=1}^ND_k$ and the contour
$\Gamma=\cup_{k=1}^N\gamma_k$.

Define {\it the Neumann jump operator} $R:C^\infty(\Gamma)\rightarrow
C^\infty(\Gamma)$ by
$$R(f)|_{\gamma_k}=\partial_\nu (V^{-}_k-V^{+}_k),$$
where $\nu$ is the outward normal to $\gamma_k=\partial D_k$, the functions
$V^{-}_k$ and $V^{+}$ are the solutions of the boundary value problems
$\Delta^{\g} V^{-}_k=0$ in $D_k$, $V^{-}|_{\partial D_k}=f$ and $\Delta^{\g}
V^{+}=0$ in $\Sigma$, $V^{+}|_{\Gamma}=f$. The Neumann jump operator is an
elliptic pseudodifferential operator of order $1$, and it is known that one
can define its determinant via the standard $\zeta$-regularization.

In what follows it is crucial that the Neumann jump operator does not change
if we vary the metric  within the same conformal class.

 Let
$(\Delta^{\g}|D_k)$ and $(\Delta^\g|\Sigma)$ be the operators of the
Dirichlet boundary problem for $\Delta^{\g}$ in domains $D_k$ and $\Sigma$
respectively, the determinants of these operators also can be defined via
$\zeta$-regularization.

 Due to Theorem
$B^{*}$ from \cite{BFK}, we have
\begin{equation}\label{s1}
{\rm  det}\Delta^{\g}=\left\{\prod_{k=1}^N{\rm
det}(\Delta^\g|D_k)\right\}\,{\rm det}(\Delta^\g|\Sigma)\,{\rm det}R\,\{{\rm
Area}(\L,\g)\}\,\{l(\Gamma)\}^{-1},
\end{equation}
where $l(\Gamma)$ is the length of the contour $\Gamma$ in the metric $\g$
\begin{remark}{\rm We have excluded the zero modes of an operator from the definition of
its determinant, so we are using  the same notation ${\rm det}\,A$ for the
determinants of operators $A$ with and without zero modes. In \cite{BFK} the
determinant of an operator $A$ with zero modes is always equal to zero, and
what we  call here ${\rm det}\, A$  is called the modified
determinant in \cite{BFK} and denoted there by ${\rm det}^* \,A$.  }.\end{remark}

An analogous statement holds for the flat conical metric. Namely let  $\L$  be a
compact polyhedral surface with vertices $P_1, \dots, P_N$ and $g$ be a
corresponding flat metric with conical singularities.  Choose the domains
$D_k$, $k=1, \dots, N$ being (open) nonoverlapping disks centered at $P_k$
and let $(\Delta|D_k)$ be the Friedrichs extension of the Laplacian with
domain $C^\infty_0(D_k\setminus P_k)$ in $L_2(D_k)$. Then formula (\ref{s1})
is still valid with $\Delta^\g=\Delta$ (cf. \cite{Leipzig1} or the recent
paper \cite{Loya} for a more general result).
\subsection{Polyakov's formula}

 We state this result in the form given in (\cite{Fay92}, p. 62).
 Let $\g_1=\rho_1^{-2}(z, \bar z)\widehat{dz}$ and $\g_2=\rho_2^{-2}(z, \bar z)\widehat{dz}$ be two
{\it smooth} conformal metrics on $\L$ and let ${\rm det}\Delta^{\g_1}$ and
${\rm det}\Delta^{\g_2}$ be the determinants of the corresponding Laplacians
(defined via the standard Ray-Singer regularization). Then
\begin{equation}\label{Polyakov}
\frac{{\rm det}\Delta^{\g_2}}{{\rm det}\Delta^{\g_1}}=\frac{{\rm Area}(\L,
\g_2)}{{\rm Area}(\L, \g_1)}
\exp\left\{\frac{1}{3\pi}\int_\L\log\frac{\rho_2}{\rho_1}\partial^2_{z\bar
z}\log(\rho_2\rho_1)\widehat{dz}\right\}\, .
\end{equation}
\subsection{Analog of Polyakov's formula for a pair of flat conical
metrics}\index{Polykov formula for flat singular metrics}
\begin{proposition}\label{poland}
Let $a_1, \dots, a_N$ and $b_1, \dots, b_M$ be real numbers which are greater
than $-1$ and satisfy $a_1+\dots+a_N=b_1+\dots+b_M=2g-2$.
 Let also $T$ be a connected $C^1$-manifold and let
$$T\ni t\mapsto \g_1(t),\ \ T\ni t\mapsto \g_2(t)$$
be two $C^1$-families of flat conical metrics on $\L$ such that
\begin{enumerate}
\item For any $t\in T$ the metrics $\g_1(t)$ and $\g_2(t)$ define the same
conformal structure on $\L$,
\item ${\g_1}(t)$ has conical singularities at $P_1(t), \dots, P_N(t)\in\L$
with conical angles $2\pi(a_1+1)$, $\dots$, $2\pi(a_N+1)$\,.
\item ${\g_2}(t)$ has conical singularities at $Q_1(t), \dots, Q_M(t)\in L$
with conical angles $2\pi(b_1+1)$, $\dots$, $2\pi(b_M+1)$\,,
\item For any $t\in T$ the sets $\{P_1(t), \dots, P_N(t)\}$ and $\{Q_1(t),
\dots, Q_M(t)\}$ do not intersect.
\end{enumerate}
Let $x_k$ be distinguished local parameter for ${\g_1}$ near $P_k$ and $y_l$
be distinguished local parameter for ${\g_2}$ near $Q_l$ (we omit the
argument $t$).

 Introduce
the functions $f_k$, $g_l$  and the complex numbers ${\bf f_k}$, ${\bf g_l}$
by
$${\g_2}=|f_k(x_k)|^2|dx_k|^2\ \  \mbox{near} \ \ P_k;\ \ \ \ \ {\bf
f_k}:=f_k(0),$$
$${\g_1}=|g_l(y_l)|^2|dy_l|^2\ \  \mbox{near} \ \  Q_l;\ \ \ \ \ {\bf
g_l}:=g_l(0).$$ Then the following equality holds

\begin{equation}\label{ConPol}\frac{{\rm det}\Delta^{\g_1}}{{\rm det}\Delta^{\g_2}}= C\
\frac{{\rm Area}\,(\L, {\g_1})}{{\rm Area}\,(\L, {\g_2})}\  \frac{
\prod_{l=1}^M|{\bf g_l}|^{b_l/6}}{\prod_{k=1}^N|{\bf f_k}|^{a_k/6}},
\end{equation}
where the constant $C$ is independent of $t\in T$.
\end{proposition}
{\bf Proof.} Take $\epsilon>0$ and introduce the disks $D_k(\epsilon)$, $k=1,
\dots, M+N$ centered at the points $P_1, \dots, P_N$, $Q_1, \dots, Q_M$;
$D_k(\epsilon)=\{|x_k|\leq \epsilon\}$ for $k=1, \dots, N$ and
$D_{N+l}=\{|y_l|\leq \epsilon\}$ for $l=1, \dots, M$.
 Let
$h_k:\overline{{\mathbb R}}_+\rightarrow {\mathbb R}$, $k=1, \dots, N+M$  be
smooth positive functions such that
\begin{enumerate}
\item
$$\int_0^1h_k^2(r)rdr=
\begin{cases}\int_0^1r^{2a_k+1}dr=\frac{1}{2a_k+2},\ \  \mbox{if} \ \  k=1, \dots, N\\
\int_0^1r^{2b_l+1}dr=\frac{1}{2b_l+2}, \ \ \mbox{if}\ \  k=N+l,\  l=1, \dots,
M
\end{cases}$$
\item
$$h_k(r)=
\begin{cases}
r^{a_k} \ \ \mbox{for}\ \  r\geq 1 \ \ \mbox{if}\ \  k=1, \dots, N\\
r^{b_l} \ \ \mbox{for} \ \ r\geq 1 \ \ \mbox{if} \ \ k=N+l, \ l=1, \dots, M
\end{cases}$$
\end{enumerate}

Define two families of {\it smooth} metrics $\g_1^\epsilon$, $\g_2^\epsilon$
on $\L$ via
$$\g_1^\epsilon(z)=\begin{cases}
\epsilon^{2a_k}h_k^2(|x_k|/\epsilon)|dx_k|^2,\ \ \ \ \ z\in D_k(\epsilon), \ \ k=1, \dots, N\\
\g(z), \ \ \ \ \ \ \ \ z\in \L\setminus\cup_{k=1}^ND_k(\epsilon)\, ,
\end{cases}$$

$$\g_2^\epsilon(z)=\begin{cases}
\epsilon^{2b_l}h_{N+l}^2(|y_l|/\epsilon)|dy_l|^2,\ \ \ \ \ z\in D_{N+l}(\epsilon), \ \ l=1, \dots, M\\
\g(z), \ \ \ \ \ \ \ \ z\in \L\setminus\cup_{l=1}^MD_{N+l}(\epsilon)\, .
\end{cases}$$

The metrics $\g_{1, 2}^\epsilon$ converge to $\g_{1, 2}$  as $\epsilon \to 0$
and $${\rm Area}(\L, \g_{1, 2}^\epsilon)={\rm Area}(\L, \g_{1, 2}).$$

\begin{lemma}\label{SOV}
 Let $\partial_t$ be the differentiation with respect to one of the
coordinates on $T$ and let ${\rm det}\Delta^{\g_{1, 2}^\epsilon}$ be the
standard $\zeta$-regularized determinant of the Laplacian corresponding to
the smooth metric $\g_{1, 2}^\epsilon$. Then
\begin{equation}\label{sovpad}
\partial_t\log{\rm  det}\Delta^{\g_{1,2}}=\partial_t\log{\rm
det}\Delta^{\g_{1, 2}^\epsilon}.
\end{equation}
\end{lemma}
To establish the lemma consider for definiteness the pair  $\g_1$ and
$\g_1(\epsilon)$. Due to the analytic surgery formulas from section
\ref{AASS} one has
\begin{equation}\label{s11}
{\rm  det}\Delta^{\g_1}=\left\{\prod_{k=1}^N{\rm
det}(\Delta^{\g_1}|D_k(\epsilon))\right\}\,{\rm
det}(\Delta^{\g_1}|\Sigma)\,{\rm det}R\,\{{\rm
Area}(\L,\g_1)\}\,\{l(\Gamma)\}^{-1},
\end{equation}
\begin{equation}\label{s22}
{\rm  det}\Delta^{\g_1^\epsilon}=\left\{\prod_{k=1}^N{\rm
det}(\Delta^{\g_1^\epsilon}|D_k(\epsilon))\right\}\,{\rm
det}(\Delta^{\g_1^\epsilon}|\Sigma)\,{\rm det}R\,\{{\rm
Area}(\L,\g_1^\epsilon)\}\,\{l(\Gamma)\}^{-1},
\end{equation}
with $\Sigma=\L\setminus\cup_{k=1}^ND_k(\epsilon)$.

Notice that the variations of the logarithms of the  first factors in
the right
hand sides of (\ref{s11}) and (\ref{s22}) vanish (these factors are
independent of $t$) whereas the variations of the logarithms of all the remaining
factors coincide. This leads to (\ref{sovpad}).

By virtue of Lemma \ref{SOV} one has the relation $$\partial_t\left\{
\log\frac{{\rm det}\Delta^{\g_1}}{{\rm Area}(\L, \g_1)} - \log\frac{{\rm
det}\Delta^{\g_2}}{{\rm Area}(\L, \g_2)}\right\}=$$
\begin{equation}\label{e1}
\partial_t\left\{  \log\frac{{\rm det}\Delta^{\g_1^\epsilon}}{{\rm Area}(\L, \g_1^\epsilon)}
- \log\frac{{\rm det}\Delta^{\g_2^\epsilon}}{{\rm Area}(\L,
\g_2^\epsilon)}\right\}.
\end{equation}
By virtue of Polyakov's formula the r. h. s. of (\ref{e1}) can be rewritten
as
\begin{equation}
\sum_{k=1}^N\frac{1}{3\pi}\partial_t\int_{D_k(\epsilon)}(\log H_k)_{x_k\bar
x_k}\log|f_k|\widehat{dx_k}-$$ $$
\sum_{l=1}^M\frac{1}{3\pi}\partial_t\int_{D_{N+l}(\epsilon)}(\log
H_{N+l})_{y_l, \bar y_l}\log|g_l|\widehat{dy_l},
\end{equation}
where $H_k(x_k)=\epsilon^{-a_k}h_k^{-1}(|x_k|/\epsilon)$, $k=1, \dots, N$ and
 $H_{N+l}(y_l)=\epsilon^{-b_l}h_{N+l}^{-1}(|y_l|/\epsilon)$, $l=1, \dots,
M$. Notice that for $k=1, \dots, N$ the function $H_k$ coincides with
$|x_k|^{-a_k}$ in a vicinity of the circle $\{|x_k|=\epsilon\}$ and the Green
formula implies that
$$\int_{D_k(\epsilon)}(\log H_k)_{x_k\bar x_k}\log|f_k|\widehat{dw_k}=\frac{i}{2}\left\{
\oint_{|x_k|=\epsilon}(\log|x_k|^{-a_k})_{\bar x_k}\log|f_k|d\bar
x_k+\right.$$
$$\left.+\oint_{|x_k|=\epsilon}\log|x_k|^{-a_k}(\log|f_k|)_{x_k}dx_k+\int_{D_k(\epsilon)}(\log|f_k|)_{x_k\bar x_k}
\log H_kdx_k\wedge d\bar x_k \right\}$$ and, therefore,
\begin{equation}\label{e2}\partial_t\int_{D_k(\epsilon)}(\log H_k)_{x_k\bar x_k}\log|f_k|\widehat{dx_k}
=-\frac{a_k\pi}{2}\partial_t\log|{\bf f}_k|+o(1)\end{equation} as
$\epsilon\to 0$. Analogously
\begin{equation}\label{e22}\partial_t\int_{D_{N+l}(\epsilon)}(\log H_{N+l})_{y_l\bar y_l}\log|g_l|\widehat{dy_l}
=-\frac{b_l\pi}{2}\partial_t\log|{\bf g}_l|+o(1)\end{equation} as
$\epsilon\to 0$.

Formula (\ref{ConPol}) follows from (\ref{e1}), (\ref{e2}) and (\ref{e22}).
$\square$

\subsection{Lemma on three polyhedra}

For any metric $\g$ on $\L$ denote by $Q(\g)$ the ratio ${\det \Delta^\g}/{\rm Area}(\L, \g)$.

 Consider three families of flat conical metrics ${\bf l}(t)\sim{\bf m}(t)\sim{\bf n}(t)$  on $\L$ (here $\sim$ means conformal equivalence), where the metric ${\bf l}(t)$ has conical points $P_1(t), \dots, P_L(t)$ with conical angles $2\pi(a_1+1), \dots, 2\pi(a_L+1)$,
the metric ${\bf m}(t)$ has conical points $Q_1(t), \dots, Q_M(t)$ with conical angles $2\pi(b_1+1), \dots, 2\pi(b_M+1)$ and
the metric ${\bf n}(t)$ has conical points $R_1(t), \dots, R_N(t)$ with conical angles $2\pi(c_1+1), \dots, 2\pi(c_N+1)$.

Let $x_k$ be the distinguished local parameter for ${\bf l}(t)$ near $P_k(t)$ and let ${\bf m}(t)=|f_k(x_k)|^2|dx_k|^2$ and ${\bf n}(t)=|g_k(x_k)|^2|dx_k|^2$
near $P_k(t)$. Let $\xi$ be an arbitrary conformal local coordinate in a vicinity of the point $P_k(t)$. Then one has
${\bf m}=|f(\xi)|^2|d\xi|^2$ and ${\bf n}=|g(\xi)|^2|d\xi|^2$ with some holomorphic functions $f$ and $g$ and the
ratio $$\frac{ {\bf m}(t)} {{\bf n}(t)}\left(P_k(t)\right):=\frac{|f(0)|^2}{|g(0)|^2}$$
is independent of the choice of the conformal local coordinate. In particular it coincides with the ratio $|f_k(0)|^2/|g_k(0)|^2$.

From Proposition \ref{poland},  one gets
the relation
$$1=\left\{\frac{Q({\bf l}(t))}{Q({\bf m}(t))}\frac{Q({\bf m}(t))}{Q({\bf n}(t))}\frac{Q({\bf n}(t))}{Q({\bf l}(t))}\right\}^{-12}=$$
\begin{equation}\label{TRI}{\rm C}\,
\prod_{i=1}^{N}\left[\frac{{\bf l}(t)}{{\bf m}(t)}(R_i(t))\right]^{c_i}
\prod_{j=1}^{L}\left[\frac{{\bf m}(t)}{{\bf n}(t)}(P_j(t))\right]^{a_j}
\prod_{k=1}^{M}\left[\frac{{\bf n}(t)}{{\bf l}(t)}(Q_k(t))\right]^{b_k}\,,\end{equation}
where the constant $C$ is independent of $t$.

From the following statement (which we call {\it the lemma on three polyhedra}) one can see that the constant $C$ in (\ref{TRI})  is equal to $1$.
\begin{lemma}\label{Ph3} Let $\L$ be a compact Riemann surface of an arbitrary genus $g$ and let
${\bf l}$, ${\bf m}$ and ${\bf n}$  be three conformal flat conical metrics on $\L$. Suppose that the metric ${\bf l}$ has conical points $P_1, \dots, P_L$ with conical angles $2\pi(a_1+1), \dots, 2\pi(a_L+1)$,
the metric ${\bf m}$ has conical points $Q_1, \dots, Q_M$ with conical angles $2\pi(b_1+1), \dots, 2\pi(b_M+1)$ and
the metric ${\bf n}$ has conical points $R_1, \dots, R_N$ with conical angles $2\pi(c_1+1), \dots, 2\pi(c_N+1)$. (All the points $P_l$, $Q_m$, $R_n$ are supposed to be distinct.) Then one has the relation
\begin{equation}\label{Poly3}
\prod_{i=1}^{N}\left[\frac{{\bf l}}{{\bf m}}(R_i)\right]^{c_i}
\prod_{j=1}^{L}\left[\frac{{\bf m}}{{\bf n}}(P_j)\right]^{a_j}
\prod_{k=1}^{M}\left[\frac{{\bf n}}{{\bf l}}(Q_k)\right]^{b_k}=1\,.
\end{equation}
\end{lemma}

{\bf Proof.} When $g>0$ and all three metrics ${\bf l}$, ${\bf m}$ and ${\bf n}$  have trivial holonomy, i. e. one has ${\bf l}=|\omega_1|^2$, ${\bf m}=|\omega_2|^2$ and ${\bf n}=|\omega_3|^2$ with some holomorphic one-forms $\omega_1$, $\omega_2$ and $\omega_3$, relation (\ref{Poly3}) is an immediate consequence of the Weil reciprocity law (see \cite{GriffHar}, \S 2.3). In general case the statement reduces to an analog of the Weil reciprocity law for harmonic functions with isolated singularities.

\section{Polyhedral tori}
Here we establish a formula for the determinant of the Laplacian on a
polyhedral torus, i.e., a Riemann surface of genus one with flat conical
metric. We do this by comparing this determinant with the determinant of the
Laplacian corresponding to the {smooth} flat metric on the same torus. For
the latter Laplacian  the spectrum is easy to find and the determinant is
explicitly known (it is given by the Ray-Singer formula stated below).

In this section $\L$ is an elliptic ($g=1$) curve and it is assumed that $\L$
is the quotient of the complex plane ${\mathbb C}$ by the lattice generated
by $1$ and $\sigma$, where $\Im\sigma>0$. The differential $dz$ on ${\mathbb
C}$ gives rise to a holomorphic differential $v_0$ on $\L$ with periods $1$
and $\sigma$.

\subsubsection{Ray-Singer formula} Let $\Delta$ be the
Laplacian on $\L$ corresponding to the flat smooth metric $|v_0|^2$.
 The
following formula for ${\rm det}\Delta$ was proved in \cite{Ray}:
\begin{equation}\label{RS}
{\rm det}\Delta=C|\Im \sigma|^2|\eta(\sigma)|^4,
\end{equation}
where $C$ is a $\sigma$-independent constant and $\eta$ is the Dedekind
eta-function.
\subsection{Determinant of the Laplacian on a polyhedral torus}

Let $\sum_{k=1}^Nb_kP_k$ be a generalized divisor on $\L$ with
$\sum_{k=1}^Nb_k=0$ and assume that $b_k>-1$ for all $k$. Let $\g$ be a flat
conical metric corresponding to this divisor via Troyanov's theorem. Clearly,
it has a finite area and is defined uniquely when this area is fixed. Fixing
numbers $b_1, \dots, b_N>-1$ such that $\sum_{k=1}^Nb_k=0$,  we define the
space ${\cal M}(b_1, \dots, b_N)$ as the moduli space of pairs $(\L, \g)$,
where $\L$ is an elliptic curve and $\g$ is a flat conformal metric on $\L$
having $N$ conical singularities with conical angles $2\pi (b_k+1)$, $k=1,
\dots, N$. The space ${\cal M}(b_1, \dots, b_N)$ is a connected orbifold of
real dimension $2N+3$.

We are going to give an explicit formula for the function
$${\cal M}(\beta_1, \dots, \beta_N)\ni (\L, \g)\mapsto {\rm det}\Delta^\g\, .$$

Write the  normalized holomorphic differential $v_0$ on the elliptic curve
$\L$ in the distinguished local parameter $x_k$ near the conical point $P_k$
($k=1, \dots, N$) as
$$v_0=f_k(x_k)dx_k$$
and define
\begin{equation}{\bf f}_k:=f_k(x_k)|_{x_k=0}, \ k=1, \dots, N\,.   \end{equation}

 \begin{theorem}\label{MT}
The following formula holds true
\begin{equation}\label{fMT}
{\rm det}\Delta^\g=C|\Im \sigma|\,{\rm Area}(\L, \g)\,
|\eta(\sigma)|^4\prod_{k=1}^N|{\bf f}_k|^{-b_k/6},
\end{equation}
where $C$ is a constant depending only on $b_1, \dots, b_N$.
\end{theorem}

{\bf Proof.} The theorem immediately follows from (\ref{RS}) and
(\ref{ConPol}).

\section{Polyhedral surfaces of higher genus}

Here we generalize the results of the previous section to the case of
polyhedral surfaces of an arbitrary genus.
 Among all polyhedral surfaces of
genus $g\geq 1$ we distinguish  {\it flat surfaces with trivial holonomy}. In
our calculation of the determinant of the Laplacian, it is this class of
surfaces which plays the role of the smooth flat tori in genus one. For
flat surfaces with trivial holonomy we find an explicit expression for the
determinant of the Laplacian which generalizes the Ray-Singer formula
(\ref{RS}) for smooth flat tori. As we did in genus one, comparing two
determinants of the Laplacians by means of Proposition \ref{poland}, we
derive a formula for the determinant of the Laplacian on a general polyhedral
surface.
\subsection{Flat surfaces with trivial holonomy and moduli spaces of holomorphic differentials on Riemann surfaces}
\index{flat surfaces with trivial holonomy} We follow \cite{ZK} and Zorich's
survey \cite{Zorich}. Outside the vertices a Euclidean polyhedral surface
$\L$ is locally isometric to a Euclidean plane and one can define the
parallel transport along paths on the punctured surface $\L\setminus\{P_1,
\dots, P_N\}$. The parallel transport along a homotopically nontrivial loop
in $\L\setminus\{P_1, \dots, P_N\}$ is generally nontrivial. If, e.g., a
small loop encircles a conical point $P_k$ with conical angle
$\beta_k$, then
a tangent vector to $\L$ turns by $\beta_k$ after the parallel transport
along this loop.

A  Euclidean polyhedral surface $\L$ is called {\it a surface with trivial
holonomy} if the parallel transport along any loop in $\L\setminus\{ P_1,
\dots, P_N\}$ does not change tangent vectors to $\L$ .

All conical points of a surface with trivial holonomy must have conical
angles which are integer multiples of $2\pi$.

A flat conical metric $g$ on a compact real oriented two-dimensional manifold
$\L$ equips $\L$ with the structure of a compact Riemann surface, if this
metric has trivial holonomy then it necessarily has the form $g=|w|^2$, where
$w$ is a holomorphic differential on the Riemann surface $\L$ (see
\cite{Zorich}). The holomorphic differential $w$ has zeros at the conical
points of the metric $g$. The multiplicity of the zero at the point $P_m$
with the conical angle $2\pi (k_m+1)$ is equal to $k_m$ \footnote{ There
exist polyhedral surfaces with nontrivial holonomy whose conical angles are
all integer multiples of $2\pi$. To construct an example take a compact
Riemann surface $\L$ of genus $g>1$ and choose $2g-2$ points $P_1, \dots,
P_{2g-2}$ on $\L$ in such a way that the divisor $P_1+\dots+P_{2g-2}$ is not
in the canonical class. Consider the flat conical
conformal metric $\g$ corresponding to the divisor $P_1+\dots + P_{2g-2}$
according to the Troyanov theorem. This metric must have nontrivial holonomy
and all its conical angles are equal to $4\pi$.}.

The holomorphic differential $w$ is defined up to a unitary complex factor.
This ambiguity can be avoided if the surface $\L$ is provided with a
distinguished direction (see \cite{Zorich}), and it is assumed that $w$ is
real along this distinguished direction. In what follows we always assume
that surfaces with trivial holonomy are provided with such a direction.

Thus, to a Euclidean polyhedral surface of genus $g$ with trivial holonomy we
put into correspondence a pair $(\L, w)$, where $\L$ is a compact Riemann
surface and $\omega$ is a holomorphic differential on this surface. This
means that we get an element  of {\it the moduli space}, ${\cal H}_g$, {\it
of holomorphic differentials over Riemann surfaces of genus} $g$ (see
\cite{ZK}).

 The space ${\cal H}_g$
is stratified according to the multiplicities of zeros of $w$.

Denote by ${\cal H}_{g}(k_1,\dots, k_M)$ the stratum of ${\cal H}_g$,
consisting of differentials  $w$ which have $M$ zeros on $\L$ of
multiplicities $(k_1,\dots,k_M)$. Denote the zeros of $w$ by $P_1,\dots,P_M$;
then the divisor  of the differential $w$ is given by $(w)=\sum_{m=1}^M k_m P_m$.
Let us  choose a  canonical basis of cycles $(a_{\a},b_{\a})$ on the Riemann
surface $\L$ and cut $\L$ along these cycles starting at the same point  to
get the fundamental polygon $\hat{\L}$. Inside $\hat{\L}$ we choose $M-1$
( homology classes of) paths
 $l_{m}$ on $\L\setminus (w)$
connecting the zero $P_1$ with other zeros $P_m$ of $w$, $m=2,\dots, M$. Then
the local coordinates on ${\cal H}_{g}(k_1,\dots, k_M)$ can be chosen as
follows \cite{KZ1}:

\be A_\a:=\oint_{a_\a} w\;,\ \  B_\a:=\oint_{b_\a} w\;,\ \  z_{m} :=\int_{l_{m}} w\;,\ \
\a=1, \dots, g;\ m=2,\dots,M \;. \la{coordint} \ee

The area of the surface $\L$ in the metric $|w|^2$ can be expressed
 in terms of these coordinates as follows:
$${\rm Area}(\L, |w|^2) =\Im \sum_{\a=1}^g A_\a \bar{B_\a}\;.$$
If all zeros of $w$ are simple, we have $M=2g-2$; therefore, the dimension of
the highest stratum ${\cal H}_{g}(1,\dots, 1)$ equals $4g-3$.

 The Abelian
integral $z(P)=\int_{P_1}^P w$ provides a local coordinate in a neighborhood
of any point $P\in \L$ except the zeros $P_1,\dots,P_M$. In a neighborhood of
$P_m$ the local coordinate can be chosen to be $(z(P)-z_m)^{1/(k_m+1)}$.

\begin{remark}{\rm The following construction  helps to visualize these coordinates in the
case of the highest stratum $H_g(1, \dots, 1)$.

Consider  $g$ parallelograms $\Pi_1, \dots, \Pi_g$ in the complex plane with
coordinate $z$ having the sides $(A_1, B_1)$, $\dots$, $(A_g, B_g)$.  Provide
these parallelograms with a system of cuts
$$[0, z_2],\ \ \  [z_3, z_4],\ \ \  \dots,\ \ \  [z_{2g-3}, z_{2g-2}]$$
(each cut should be repeated on two different parallelograms). Identifying
opposite sides of the parallelograms and gluing the obtained $g$ tori along
the cuts, we get a compact Riemann surface $\L$ of genus $g$. Moreover, the
differential $dz$ on the complex plane gives rise to a holomorphic
differential $w$ on $\L$ which has $2g-2$ zeros at the ends of the cuts.
Thus, we get a point $(\L, w)$ from ${\cal H}_g(1, \dots, 1)$. It can be
shown that any generic point of ${\cal H}_g(1, \dots, 1)$ can be obtained via
this construction; more sophisticated gluing is required to represent points
of other strata, or non generic points of the stratum ${\cal H}_g(1, \dots,
1)$. }\end{remark} To shorten the notations it is convenient to consider the
coordinates $A_\a$,$B_\a$, $z_m$ altogether. Namely,  in the sequel we
shall denote them by $\ko_k$, $k=1,\dots,2g+M-1$, where \be\label{sokr}
\ko_\a:= A_\a\;,\  \ \ko_{g+\a}:= B_\a\;,\ \  \a=1,\dots, g\;,\ \  \ko_{2g+m}:=
z_{m+1}\;,\  m= 1,\dots, M-1 \ee

Let us also introduce corresponding cycles $\cy_k$, $k=1,\dots,2g+M-1$, as
follows: \be \cy_\a= - b_\a\;,\ \  \cy_{g+\a}= a_\a\;,\ \  \a=1,\dots, g\;; \ee the
cycle $\cy_{2g+m}$, $m= 1,\dots, M-1$ is defined to be the small circle with
positive orientation around the point $P_{m+1}$.

\subsubsection{Variational formulas on the spaces of holomorphic
differentials} In the previous section we introduced the coordinates on the
space of surfaces with trivial holonomy and fixed type of conical
singularities. Here we study the behavior of basic objects on these surfaces
under the change of the coordinates. In particular, we derive  variational
formulas of Rauch type for the matrix of ${\bf b}$-periods of the
underlying Riemann surfaces. We also give variational formulas for the Green
function, individual eigenvalues, and the determinant of the Laplacian on
these surfaces.

{\bf Rauch formulas on the spaces of holomorphic differentials.}\index{Rauch
formulas} \ For any compact Riemann surface $\L$ we introduce the prime-form
$E(P,Q)$ and the canonical meromorphic bidifferential
\begin{equation}\label{bergdef}
\Be(P,Q)= d_P d_Q\log E(P,Q)\end{equation} (see \cite{Fay92}).
 The bidifferential $\Be(P,Q)$ has
the following local behavior as $P\to Q$: \be \Be(P,Q)=
\left(\f{1}{(x(P)-x(Q))^2}+\f{1}{6} S_B(x(P))+ o(1)\right)dx(P)dx(Q),
\la{defproco} \ee where $x(P)$ is a local parameter. The term $S_B(x(P))$ is
a projective connection which is called {\it the Bergman projective
connection} (see \cite{Fay92}).

Denote by $\hd_\a(P)$ the basis of holomorphic 1-forms on $\L$ normalized by
$\int_{a_\a} \hd_\b=\delta_{\a\b}$.

The matrix of {\bf b}-periods of the surface $\L$ is given by
$\B_{\a\b}:=\oint_{b_\a} \hd_\b$.

\begin{proposition}\la{varfow}(see \cite{Leipzig}) Let a pair $(\L, w)$ belong to the space ${\cal H}_{g}(k_1,\dots,
k_M)$. Under variations of the coordinates on ${\cal H}_{g}(k_1,\dots, k_M)$
the normalized holomorphic differentials and the matrix of ${\bf b}$-periods
of the surface $\L$ behaves as follows:
 \be \f{\p \hd_\a(P)}{\p
\ko_k}\Big|_{z(P)} = \f{1}{2\pi i}\oint_{\cy_k}\f{\hd_\a(Q)\Be(P,Q)}{w(Q)}\;,
\la{varw} \ee \be\label{varB} \f{\p\B_{\a\b}}{\p
\ko_k}=\oint_{\cy_k}\f{\hd_\a \hd_\b}{w}\ee
 where $k=1,\dots,
2g+M-1$;  we assume that the local coordinate
 $z(P)=\int_{P_1}^{P}w$
is kept constant under differentiation.
\end{proposition}

{\bf Variation of the resolvent kernel and eigenvalues.}\ For a pair $(\L,
w)$ from ${\cal H}_{g}(k_1,\dots, k_M)$ introduce the Laplacian
$\Delta:=\Delta^{|w|^2}$ in the flat conical metric $|w|^2$ on $\L$ (recall that
we always deal with the Friedrichs extensions). The corresponding resolvent
kernel $G(P, Q; \l)$, $\l\in {\mathbb C}\setminus {\rm sp}\,(\Delta)$
\begin{itemize}
\item satisfies $(\Delta_P-\l)G(P, Q; \l)=(\Delta_Q-\l)G(P, Q; \l)=0$ outside the diagonal
$\{P=Q\}$,
\item is bounded near the conical points i. e. for any $P\in \L\setminus
\{P_1, \dots, P_M\}$
$$G(P, Q; \l)=O(1)$$
as $Q\to P_k$, $k=1, \dots, M$,
\item obeys the asymptotics $$G(P,
Q; \l)=\frac{1}{2\pi}\log|x(P)-x(Q)|+O(1)$$ as $P\to Q$, where $x(\cdot)$ is
an arbitrary (holomorphic) local parameter near $P$.
\end{itemize}
The following proposition is an analog of the classical Hadamard formula for
the variation of the Green function of the Dirichlet problem in a plane
domain.\index{resolvent kernel}

\begin{proposition} The following variational formulas for the resolvent
kernel $G(P, Q; \l)$ hold:
\begin{equation}\label{Had1}
\frac{\partial G(P, Q; \l)}{\partial A_\alpha}=2i\int_{{\bf b}_\alpha}\omega(P, Q;
\lambda)\,,\end{equation}

\begin{equation}\label{Had2}
\frac{\partial G(P, Q; \l)}{\partial B_\alpha}=-2i\int_{{\bf a}_\alpha}\omega(P, Q;
\lambda)\,,\end{equation} where
$$\omega(P, Q; \l)=G(P, z; \l)G_{z\bar{z}}(Q, z; \l)\overline{dz}+G_z(P, z;
\l)G_z(Q, z; \l)dz$$ is a closed $1$-form and $\alpha=1, \dots, g$;

\begin{equation}\label{Had3}
\frac{\partial G(P, Q; \l)}{\partial z_m}=-2i\lim_{\epsilon\to
0}\oint_{|z-z_m|=\epsilon}G_z(z, P; \l)G_z(z, Q; \l)dz\,,\end{equation} where
$m=2, \dots, M$. It is assumed that the coordinates $z(P)$ and $z(Q)$ are
kept constant under variation of the moduli $A_\alpha, B_\alpha, z_m$.
\end{proposition}
\begin{remark}{\rm One can unite the formulas (\ref{Had1}-\ref{Had3}) in a
single formula:
$$\frac{\partial G(P, Q; \l)}{\partial \ko_k}=$$
\begin{equation}\label{Had}
-2i\left\{\int_{\cy_k}\frac{G(R, P; \l)\partial_R \overline{\partial_R} G(R,Q;
\l)+\partial_RG(R, P; \l)\partial_R G(R,Q; \l)}{w(R)}\right\}\,,
\end{equation}
where $k$=1, \dots, 2g+M-1. }
\end{remark}
{\bf Proof.}
 We start with the following integral representation of a solution $u$
 to the homogeneous equation $\Delta u-\l u=0$ inside the fundamental polygon $\hat{\L}$:
\begin{equation}\label{green}
u(\xi, \bar \xi)=-2i\int_{\partial \hat{\L}}G(z, \bar z, \xi, \bar\xi;
\l)u_{\bar z}(z, \bar z)d\bar z+G_z(z, \bar z, \xi, \bar \xi; \l)u(z, \bar
z)dz\,.
\end{equation}
Cutting the surface $\L$ along the basic cycles, we notice that the function
$\dot{G}(P,\  \cdot\ ; \l)=\frac{\partial G(P,\  \cdot\ ; \l)}{\partial
B_\beta}$ is a solution to the homogeneous equation $\Delta u-\l u=0$ inside
the fundamental polygon (the singularity of $G(P, Q; \l)$ at $Q=P$ disappears
after differentiation) and that the functions $\dot{G}(P,\ \cdot\ ; \l)$ and
$\dot{G}_{\bar{z}}(P,\  \cdot\ ; \l)$ have the jumps $G_z(P,\ \cdot\ ; \l)$
and $G_{z\bar{z}}(P, \ \cdot\ ; \l)$ on the cycle ${\bf a}_\beta$. Applying
(\ref{green}) with $u=\dot{G}(P,\  \cdot\ ; \l)$, we get (\ref{Had2}).
Formula (\ref{Had1}) can be proved in the same manner.

The relation $d\omega(P, Q; \lambda)=0$ immediately follows from the equality
$G_{z\bar z}(z, \bar z, P; \l)=\frac{\l}{4} G(z, \bar z, P; \l)$.

Let us prove (\ref{Had3}). From now on we assume for simplicity that $k_m=1$,
where $k_m$ is the multiplicity of the zero $P_m$ of the holomorphic
differential $w$.

Applying  Green's formula (\ref{green}) to the domain
$\hat{\L}\setminus\{|z-z_m|<\epsilon\}$ and $u=\dot{G}=\frac{\partial
G}{\partial z_m}$, one gets $$\dot{G}(P, Q; \l)=$$
\begin{equation}\label{tvetv}
2i\lim_{\epsilon\to 0}\oint_{|z-z_m|=\epsilon}\dot{G}_{\bar z}(z, \bar z, Q;
\l)G(z, \bar z, P; \l)\bar {dz}+\dot{G}(z, \bar z, Q; \l)G_z(z, \bar z, P;
\l)dz\,.
\end{equation}

Observe that the function $x_m \mapsto G(x_m, \bar x_m, P; \l)$ (defined in a
small neighborhood of the point $x_m=0$) is a bounded solution to the
elliptic equation $$\frac{\partial^2 G(x_m, \bar x_m, P; \l)}{\partial
x_m\partial \bar{x}_m} -\lambda|x_m|^2G(x_m, \bar x_m, P; \l)=0$$ with real
analytic coefficients and, therefore, is real analytic near $x_m=0$.

From now on we write $x$ instead of $x_m=\sqrt{z-z_m}$. Differentiating the
expansion
\begin{equation}\label{expan}
G(x, \bar x, P; \l)=a_0(P, \l)+a_1(P, \l)x+a_2(P, \l)\bar x+a_3(P, \l)x\bar
x+\dots
\end{equation}
with respect to $z_m$, $z$ and $\bar z$, one gets the asymptotics
\begin{equation}\label{assy11}\dot{G}(z, \bar z, Q; \l)=-\frac{a_1(Q,
\l)}{2x}+O(1),\end{equation} \begin{equation}\label{assy2}\dot{G}_{\bar z}(z,
\bar z, Q; \l)=\frac{\dot{a}_2(Q, \l)}{2\bar x}-\frac{a_3(Q, \l)}{4x\bar
x}+O(1),\end{equation} \begin{equation}\label{assy3}G_z(z, \bar z, P;
\l)=\frac{a_1(P, \l)}{2x}+O(1),\end{equation} Substituting (\ref{assy11}),
(\ref{assy2}) and (\ref{assy3}) into (\ref{tvetv}), we get the relation
$$\dot{G}(P, Q, \l)=2\pi a_1(P, \l)a_1(Q, \l).$$
On the other hand, calculation of the right hand side of formula (\ref{Had3})
via (\ref{assy3}) leads to the same result. $\square$

Now we give a variation formula for an eigenvalue of the Laplacian on a flat
surface with trivial holonomy.

\begin{proposition}

Let $\l$ be an eigenvalue of $\Delta$ (for simplicity  we assume it to have
multiplicity one) and let $\phi$ be the corresponding normalized
eigenfunction. Then
\begin{equation}\label{sobstv}
\frac{\partial \lambda}{\partial \ko_k}=2i\int_{\cy_k}\left(\frac{(\partial
\phi)^2}{w}+\frac{1}{4}\lambda\phi^2\bar w\right)\,,
\end{equation}
where $k=1,\dots, 2g+M-1$.
\end{proposition}
{\bf Proof.}  For brevity we give the proof only for the case $k=g+1, \dots,
2g$. One has
$$\iint_{\hat{L}}\phi\dot{\phi}=\frac{1}{\l}\iint_{\hat{\L}}\Delta
\phi\,\dot{\phi}=\frac{1}{\l}\left\{2i\int_{\partial \hat{\L}}(\phi_{\bar
z}\dot{\phi}d\bar z+\phi\dot{\phi}_z\,dz)+\iint_{\hat{\L}}\phi(\l\phi)^\cdot
\right\}=$$ $$\frac{1}{\lambda}\left\{ 2i\int_{{\bf a}_\beta}(\phi_{\bar
z}\phi_z\,d\bar
z+\phi\phi_{zz}\,dz)+\dot{\l}+\l\iint_{\hat{\L}}\phi\dot{\phi}
 \right\}\,.$$
This implies (\ref{sobstv}) after integration by parts (one has to make use
of the relation $d(\phi \phi_z)=\phi_z^2dz+\phi\phi_{zz}dz+\phi_{\bar
z}\phi_zd\bar z+\frac{1}{4}\l \phi^2d\bar z$).  $\square$

{\bf Variation of the determinant of the Laplacian.} For simplicity we
consider only flat surfaces with trivial holonomy having $2g-2$ conical
points with conical angles $4\pi$. The proof of the following proposition can
be found in  \cite{Leipzig}.

\begin{proposition}\label{p1} Let $(\L, w)\in {\cal H}_g(1, \dots, 1)$.
Introduce the notation \be {\mathbb Q}(\L, |w|^2):= \Big\{\frac{{\rm  {det}}
\, \Delta^{|w|^2}}{{\rm Area}(\L, |w|^2)\,{\rm det}\Im \B} \Big\}\; \la{TLW0}
\ee
 where  $\B$ is the matrix of ${\bf b}$-periods of the surface $\L$ and
${\rm Area}(\L,|w|^2)$ denotes the area of  $\L$ in the metric $|w|^2$.

The following variational formulas hold
\begin{equation}\label{cor1}
\frac{\partial \log {\mathbb Q}(\L, |w|^2)}{\partial \ko_k}=-\frac{1}{12\pi
i}\oint_{\cy_k}\frac{S_B-S_w}{w}\;,
\end{equation}
where $k=1,\dots, 4g-3$;
 $S_B$ is the Bergman projective connection,
$S_w$ is the projective connection given by the Schwarzian derivative
$\Big\{\int^Pw, x(P) \Big\}$; $S_B-S_w$ is a meromorphic quadratic
differential with poles of the second order at the zeroes $P_m$ of $w$.
\end{proposition}

\subsubsection{An explicit formula for the determinant of the Laplacian on a
flat surface with trivial holonomy} We start with recalling the properties of
the prime form $E(P, Q)$ (see \cite{Fay73,Fay92}, some of these properties
were already used in our proof of the Troyanov theorem above).

\begin{itemize}
\item The prime form $E(P, Q)$
 is an antisymmetric  $-1/2$-differential with respect to both $P$ and
$Q$,
\item Under tracing of $Q$ along the cycle ${\bf a}_\a$ the prime-form
remains invariant; under the tracing along  ${\bf b}_\alpha$ it gains the factor
\begin{equation}\label{primetwist}
\exp\left(-\pi i \B_{\a\a}-2\pi i\int_P^Q v_\a\right)\;.
\end{equation}
\item On the diagonal $Q\to P$  the prime-form has first order zero
with the following asymptotics:
$$E(x(P), x(Q))\sqrt{dx(P)}\sqrt{dx(Q)}=$$
\begin{equation}\label{primas}
(x(Q)-x(P))\left(1-\frac{1}{12}S_B(x(P))(x(Q)-x(P))^2+O((x(Q)-x(P))^3\right),
\end{equation}
where $S_B$ is the Bergman projective connection and $x(P)$ is an arbitrary
local parameter.
\end{itemize}

 The next object we shall need is the
vector of Riemann constants: \be
K^P_\alpha=\frac{1}{2}+\frac{1}{2}\B_{\a\a}-\sum_{\b=1,
\b\neq\a}^g\oint_{{\bf a}_\b}\left(\hd_\beta\int_P^x\hd_\a\right) \la{rimco}\ee
where the interior integral is taken along a path which does not intersect
$\partial\widehat\L$.

 In what follows the pivotal role is played by the following holomorphic multivalued
$g(1-g)/2$-differential on $\L$
\begin{equation}\label{c}
\cdiff(P)=\frac{1}{\Wcal[v_1, \dots, v_g](P)}\sum_{\a_1, \dots, \a_g=1}^g
\frac{\partial^g\Theta(K^P)}{\p z_{\a_1}\dots \p z_{\a_g}} v_{\a_1}\dots
v_{\a_g}(P)\;,
\end{equation}
where $\Theta$ is the theta-function of the Riemann surface $\L$, \be
\Wcal(P):= {\rm \det}_{1\leq \a, \b\leq g}||\hd_\b^{(\a-1)}(P)|| \la{Wronks}
\ee is the Wronskian determinant of holomorphic differentials
 at the point $P$.

 This differential
has multipliers $1$ and $\exp\{-\pi i (g-1)^2\B_{\a\a}-2\pi i (g-1) K_\a^P\}$
along basic
cycles ${\bf a}_\a$ and ${\bf b}_\a$, respectively.

In what follows we shall often treat tensor objects like $E(P, Q)$,
$\cdiff(P)$, etc. as scalar functions of one of the arguments (or both). This
makes sense after fixing the local system of coordinates, which is usually
taken to be $z(Q)=\int^Qw$. In particular, the expression ``the value of the
tensor $T$ at the point $Q$ in local parameter $z(Q)$"  denotes the value
of the scalar $Tw^{-\alpha}$ at the point $Q$, where $\alpha$ is the tensor
weight of $T(Q)$.

The following proposition was proved in \cite{Leipzig}.

\begin{proposition} \la{p2}
Consider the highest stratum ${\cal H}_g(1,\dots,1)$ of the space ${\cal
H}_g$ containing Abelian differentials $w$ with simple zeros.

Let us choose  the fundamental polygon $\hat{\L}$  such that
$\Abel_P((w))+2K^P=0$, where $\Abel_P$ is the Abel map with the initial point
$P$. Consider the following expression
\begin{equation}\label{otvtau}
\tau(\L,w) = {\Fcal}^{2/3}  \prod_{m,l=1\;\; m < l}^{2g-2} [E(Q_m,
Q_l)]^{{1}/{6}}\,,
\end{equation}
where the quantity \be \Fcal := [w(P)]^{\f{g-1}{2}}\c(P) \prod_{m=1}^{2g-2}[
E (P,Q_m)]^{\f{(1-g)}{2}} \ee does not depend on $P$; all prime-forms are
evaluated at the zeroes $Q_m$ of the differential $w$ in the distinguished
local parameters $x_m(P)=\left(\int_{Q_m}^P w\right)^{1/2}$. Then
\begin{equation}\label{vartau}
\frac{\partial \log \tau}{\partial \ko_k}=-\frac{1}{12\pi
i}\oint_{\cy_k}\frac{S_B-S_w}{w}\;,
\end{equation}
where $k=1, \dots, 4g-3$.
\end{proposition}

The following Theorem immediately follows from Propositions \ref{p1} and
\ref{p2}. It can be considered as a natural generalization of the Ray-Singer
formula (\ref{RS}) to the higher genus case.

\begin{theorem}\label{mmaaiinn}
Let a pair $(\L, w)$ be a point of the space ${\cal H}_g(1, \dots, 1)$. Then
the determinant of the Laplacian $\Delta^{|w|^2}$  is given by the following
expression
\begin{equation}\label{MAINDET}
{{\rm  det}}\,\Delta^{|w|^2}=C\;{\rm Area}(\L,|w|^2)\;{{\rm  det}}\Im
\B\;|\tau(\L,w)|^2,
\end{equation}
where the constant $C$  is independent of a point of ${\cal H}_g(1, \dots,
1)$. Here $\tau (\L, w)$ is given by (\ref{otvtau}).
\end{theorem}
\subsection{Determinant of the Laplacian on an arbitrary polyhedral surface of genus $g>1$}

Let $b_1, \dots, b_N$ be real numbers such that $b_k>-1$ and
$b_1+\dots+b_N=2g-2$. Denote by ${\cal M}_g(b_1, \dots, b_N)$ the moduli
space of pairs $(\L, \g)$, where $\L$ is a compact Riemann surface of genus
$g>1$ and $\g$ is a flat conformal conical metric on $\L$ having $N$ conical
points with conical angles $2\pi(b_1+1), \dots, 2\pi(b_N+1)$. The space
${\cal M}_g(b_1, \dots, b_N)$ is a (real) orbifold of (real) dimension
$6g+2N-5$.
Let $w$ be a holomorphic differential with $2g-2$ simple zeroes on $\L$.
Assume also that the set of conical points
of the metric $\g$ and the set of zeros of the differential $w$ do not
intersect.

Let $P_1, \dots, P_N$ be the conical points of $\g$ and let $Q_1, \dots,
Q_{2g-2}$ be the zeroes of $w$.
 Let $x_k$ be a distinguished local parameter for ${\g}$ near $P_k$
and $y_l$ be a distinguished local parameter for $w$ near $Q_l$.
 Introduce
the functions $f_k$, $g_l$  and the complex numbers ${\bf f_k}$, ${\bf g_l}$
by
$$|w|^2=|f_k(x_k)|^2|dx_k|^2\ \  \mbox{near} \ \ P_k;\ \ \ \ \ {\bf
f_k}:=f_k(0),$$
$$\g=|g_l(y_l)|^2|dy_l|^2\ \  \mbox{near} \ \  Q_l;\ \ \ \ \ {\bf
g_l}:=g_l(0).$$
 Then from (\ref{ConPol}) and (\ref{MAINDET}) and the lemma on three polyhedra from \S 4.4 it follows the relation
\begin{equation}\label{GenPoly}{\rm det}\Delta^{\g}=C
{\rm Area}\,(\L, {\g}) {{\rm  det}}\Im \B\;|\tau(\L,w)|^2
  \frac{ \prod_{l=1}^{2g-2}|{\bf
g_l}|^{1/6}}{\prod_{k=1}^N|{\bf f_k}|^{b_k/6}},
\end{equation}
where the constant $C$ depends only on $b_1, \dots, b_N$ (and neither the differential $w$ nor the point $(\L, \g)\in {\cal M}_g(b_1, \dots, b_N)$) and $\tau(\L, w)$ is given
by (\ref{otvtau}).

{\bf Acknowledgements.} The author is grateful to D. Korotkin for numerous suggestions, in particular, his criticism of an earlier version of this paper \cite{MPprep} lead to appearance of the lemma from \S 4.4 and a considerable improvement of our main result (\ref{GenPoly}). The author also thanks A. Zorich for very useful discussions. This paper was written during the author's stay in Max-Planck-Institut f\"ur Mathematik in Bonn, the author thanks the Institute for excellent working conditions and hospitality.


\end{document}